\documentclass[reqno]{amsart}
\usepackage{amssymb,latexsym,amscd,graphicx}
\usepackage{color}
\newtheorem{lem}{Lemma}[section]

\newtheorem{thm}{Theorem}[section]

\theoremstyle{definition}

\theoremstyle{remark}

\definecolor{orange}{rgb}{1,0.5,0}
\usepackage{graphicx}
\graphicspath{{converted_graphics/}}

\begin{document}

 \title{Automorphic Integrals with  Log-Polynomial Period Functions and Arithmetical Identities}

\author{Tewlede G/Egziabher, Hunduma Legesse Geleta and  Abdul Hassen }        
\address{Department of Mathematics: Ambo University, Ambo, Addis Ababa University, Addis Ababa, Rowan University, Glassboro, NJ 08028.}
\email{  teweldega@gmail.com, hunduma.legesse@aau.edu.et, hassen@rowan.edu }

\date{\today}          
\maketitle

\begin{abstract}  Building on the works of S. Bochner on equivalence of modular relation with functional equation associated to the Dirichlet series, K. Chandrasekharan and R. Narasimhan obtained new equivalences between the functional equation and some arithmetical identities. Sister Ann M. Heath considered the functional equation in the Hawkins and Knopp context and showed its equivalence to two arithmetical identities associated with entire modular cusp integrals involving rational period functions for the full modular group. In this paper we use techniques of Chandrasekharan and Narasimhan and extend the results of Sister Ann M. Heath to entire automorphic integrals involving rational period functions on discrete Hecke group. Moreover, we establish equivalence of two arithmetical identities with a functional equation associated with automorphic integrals involving log-polynomial-period functions on the Hecke groups. 
 
 \end{abstract}
 \maketitle
\def\theequation{\thesection.\arabic{equation}}
\section{Introduction}
\setcounter{equation}{0}

In $1951$, Bochner introduced in ~\cite{Bockner2} and  showed that the modular relation
 \[
 f(x) = x^{-\delta} g\left
 ({\frac{1}{x}}\right ) + P(x),\, 0 < x < \infty \, ,
 \] 
 where $f$ and $g$ are exponential series, $P$ is a ``residual'' function and
 $\delta  >0$, is equivalent to the functional equation
 \begin{equation}
 \label{a}
 \Gamma(s) \phi(s) = \Gamma (\delta-s)\psi(\delta-s) \, ,
 \end{equation}
 where $\Gamma(s)$ is the standard gamma function, $ \phi(s) =
 \sum^\infty_{n=1} \frac{a_n}{\lambda^s_n}$ and $\psi(s) = 
 \sum^\infty_{n=1} \frac{b_n}{\mu^s_n}$.  Here, he assumes that
 $\{\lambda_n\}$ and 
 $\{\mu_n\}$ are sequences of real numbers, 
 with
 \[
 0 < \lambda_1 < \lambda_2 < \cdots < \lambda_n \to \infty
 \]
 and 
 \[
 0 < \mu_1 < \mu_2 < \cdots < \mu_n \to \infty \, ,
 \] 
 and $\{a_n\},$ $\{b_n\}$ are two sequences of 
 complex numbers not identically zero.
 
Building on the works of Bochner,   Chandrasekharan and Narasimhan \cite{CN} showed that 
 the functional equation~(\ref{a}) is equivalent to two identities given by
 
 \begin{equation} \label{b} 
 \frac{1}{\Gamma(\rho +1)}{\sum_{\lambda_n \leq x}}^\prime
 a_n(x - \lambda_n)^\rho = \left( \frac{1}{2\pi}\right)^\rho \sum_{n=1}^\infty
 \left( \frac{x}{\mu_n} \right)^{\frac{\delta +
 \rho}{2}} b_n J _{\delta + \rho}  \{4 \pi \sqrt{\mu_n x} \} + Q_\rho (x) \, ,
 \end{equation}
 for
 $x > 0$ and $\rho \geq 2 \beta -\delta -\frac{1}{2},$ where, the prime notation $'$ on the summation sign indicates that if    $\rho=0$   and   $x=\lambda_{n}$   for some positive integer   $n$,    then we count only   $\frac{1}{2}a_{n}$,
 \[
 Q_\rho(x) =
 \frac{1}{2\pi i} \oint_c \frac{ \chi(s)(2\pi)^s x^{s+ \rho}}{\Gamma(s + \rho +
 1)}ds \, ,
 \]
 $\sum_{n=1}^\infty \frac{\vert b_n \vert}{\mu_n-\beta} < \infty$, $J_\nu(z)$ denotes the Bessel function of the first kind of order $\nu,$ and
 \begin{equation} \label{c} 
 \left( -\frac{1}{ s}\frac{d}{ds}\right)^\rho \left[ \frac{1}{s} 
    \sum^\infty_{n=1} a_n e^{-s\sqrt{\lambda_n}} \right] 
     = 2^{3\delta+\rho}
 \Gamma(\delta+\rho+\frac{1}{2})\pi^{\delta-\frac{1}{2}}
    \sum^\infty_{n=1} \frac{b_n}{ (s^2 + 16\pi^2\mu_n)^{\delta+\rho+\frac{1}{2}}}
    + R_\rho(s) \, ,
 \end{equation}
 for Re $s > 0$, where
 \[
 R_\rho(s) = \frac{1}{2\pi i} \oint_c \frac{\chi(z) (2\pi)^z \Gamma(2z +
 2\rho + 1) 2^{-\rho}} { \Gamma(z + \rho + 1)} s^{-2z-2\rho-1} dz \, ,
 \]
 for Re $s > 0,$ $\rho$ integral, $\rho  -\delta -\frac{1}{2}$ and $\rho \geq 2 \beta -\delta
 -\frac{1}{2}$.  \\
 
  Sister Ann Heath, in \cite{SISANN} showed the equivalence of the above forms for the full modular group. \\
 
  In our work here, we will show that these results hold for a more general Hecke Groups. Our notations will be different from the ones shown above. The structure of our paper is as follows: in section 2 we review some preliminary concepts  which we use in the later sections. In section 3 we state and proof our main results which are in Theorem 3.1 and Theorem 3.2. In section four we give some conclusions and discussions. 
  
\maketitle
\def\theequation{\thesection.\arabic{equation}}
 \section{ Preliminaries}\label{prelim}
\setcounter{equation}{0}

In this section, we recall some definitions, results proved elsewhere, and make some conventions for our later discussion.\\

  For    $\lambda\in\mathbb{R^{+}},$    the Hecke group    $G(\lambda)$    is defined as the subgroup of    $SL_{2}(\mathbb{R})$   given by $$G(\lambda)=\left\langle 
 \begin{pmatrix} 
 1 &\lambda \\ 0& 1
 \end{pmatrix}
 \begin{pmatrix} 0& 1 \\ -1& 0 
 \end{pmatrix}\right\rangle . $$
Note that  $G(\lambda)$ can be viewed as a group of  linear fractional transformations  generated by $S(z)=z+\lambda$    and    $T(z)=-\frac{1}{z}$. In this context, we see that $G(\lambda)$ acts on the Riemann sphere as linear fractional transformation, that is    $Mz=\frac{az+b}{cz+d}$   for    $ M=\begin{pmatrix}a & b\\ c & d \end{pmatrix}$ $ \in G(\lambda),$   and    $z\in\mathbb{C}\bigcup\lbrace \infty\rbrace,$   thus    $M$   and    $-M$   can be identified as the same matrices.  E. Hecke \cite{EHECKE}   showed the groups    $G(\lambda)$    is discrete (operate discontinuously) as linear fractional transformations on the upper half plane    $ \mathcal{H}=\{ z=x + i y: y>0\}$    if and only if either    $$\lambda>2\  \mbox{or } \ \lambda=\lambda_{p}:=2\cos\left(\frac{\pi}{p}\right),\   \mbox{with} \     3\leq p\in\mathbb{N}\cup\lbrace \infty\rbrace .$$    Clearly  $G(\lambda_3)=\Gamma(1)$    is  the  full modular group and  $G(\lambda_{\infty})=\Gamma_{\theta}$  is the theta group and is denoted by $\Gamma_{\theta}.$ \\

 A complex valued function  $\nu$  defined on  $G(\lambda)$ is called a {\em multiplier system with weight $2k$  for $G(\lambda)$} if  $|\nu(M)|=1$ for   all ~$ M= \begin{pmatrix} a & b\\ c & d \end{pmatrix} \in G(\lambda)$,and ~$\nu(M)$  satisfies  the {\em consistency condition} $$\nu(M_{3})\left(c_{3}z+d_{3} \right)^{2k}=\nu(M_{1})\nu(M_{2})\left(c_{1}M_{2}z+d_{1} \right)^{2k} \left(c_{2}z+d_{2} \right)^{2k}$$ for all $M_{1},M_{2}\in G(\lambda), M_{1}M_{2}=M_{3},$~ $ M_{j}=\begin{pmatrix}
 a_j & b_j\\c_j&d_j 
 \end{pmatrix},$
 $j=1,2,3 ~z\in\mathcal{H},$ where $k\in\mathbb{R}$.   \\
 
  A holomorphic function $F$  is called {\em entire automorphic integral of weight $2k$~with multiplier system $\nu$ for $G(\lambda)$} if it 
 admits an exponential series expansion
   \begin{equation}\label{LPP3}
 F(z)=\sum_{m=0}^{\infty} a_{m}e^{2\pi i mz/\lambda},
 \end{equation}
 where $ \Im z=y>0$~and ~$a_{m}=\mathcal{O}(m^{\delta}),$~for some ~$\delta\in\mathbb{R^{+}}$,  and satisfies the transformations
   \begin{equation}\label{LPP1}
 \bar{\nu}(S_{\lambda}) F(z+\lambda)=F(z),  ~ \nu(S_{\lambda})=e^{2\pi i\kappa}, 0\leq\kappa<1
 \end{equation}  
 
 and 
 
 \begin{equation}\label{LPP2}
 \bar{\nu}(T)z^{-2k}F\left( \frac{-1}{z} \right)=F(z)+q(z)~\forall z\in\mathcal{H},
 \end{equation}
 where ~$q(z)=\displaystyle\sum_{j=1}^{N}z^{\alpha_{j}}\sum_{t=0}^{M_{j}}\beta_{jt}\left(\log z\right)^{t},$~$\alpha_{j},~\beta_{jt}\in\mathbb{C}$.

 The function $q(z)$ in (\ref{LPP2}) is called the log-polynomial period  function of the entire automorphic integral function $F.$\\
 
 A. Hassen in \cite{ABHA}   has completely characterized  the log-polynomial period functions for entire automorphic integrals of weight $2k,~k\in\mathbb{R}$  on  the  discrete Hecke group $G(\lambda)$~for the following cases: 
 \begin{enumerate}
 	\item [(1)]
 	$k\geq 1,~\nu(S_{\lambda})=1$
 \end{enumerate}
 \begin{enumerate}
 	\item [(2)] 
 	$k>0,~\nu(S_{\lambda})\neq1$
 \end{enumerate}
 \begin{enumerate}
 	\item [(3)]
 	$k\geq0,~2k\in\mathbb{Z},~\nu(S_{\lambda})=1$~and
 \end{enumerate}
 \begin{enumerate}
 	\item [(4)]
 	$k\leq 0,~\nu(S_{\lambda})\neq 1.$
 \end{enumerate}

 For the purpose of this work we consider only the multiplier system which which $\nu(S_{\lambda})=1$ (or $(\kappa=0)$. \\

 In \cite{PAUL} Paul C.Pasles showed Riemann-Hecke-Bochner correspondence for entire automprphic integrals on the Hecke groups. We need a modified result of Pasles  \cite{PAUL} (see Theorem 3.1), which we state and proof in  
 
 \begin{thm}\label{ThM41}
 	Let $k\in\mathbb{R^{+}}\cup \lbrace 0\rbrace,~\lambda\geq2$~or~$\lambda=2\cos\left( \frac{\pi}{p}\right),~p\in\mathbb{Z},~p\geq3.$
 Suppose that $F$ is an entire automorphic integral of weight $2k,$ multiplier system $\nu$ and $\nu(S_{\lambda})=1.$ Let $q(z)$ be log-polynomial period function of $F$ on $G(\lambda)$ and let $F$ has Fourier expansion of the form $$F(z)=\sum_{m=0}^{\infty} a_{m}e^{2\pi imz/\lambda},$$ for $z\in\mathcal{H},$ 
 where
 \begin{enumerate}
 	\item[(i)] $a_{m}=\mathcal{O}(m^{\beta})$~for some ~$\beta>0,m\rightarrow\infty.$
 	\item[(ii)]Let ~$\displaystyle\varphi(s)=\sum_{m=1}^{\infty}a_{m}m^{-s},~s=\sigma+it.$~And put
 	\item[(iii)]$\displaystyle\Phi_{F}(s)=\left(\frac{2\pi}{\lambda} \right)^{-s}\Gamma(s)\varphi(s),$~for~$\sigma>\beta+1$.
 \end{enumerate}
 Then
 \begin{enumerate}
 	\item[(A)] $\Phi_{F}(s)$~has a meromorphic continuation to the whole complex plane  with at most a finite  number of poles given by:
 	
 	$$\Phi_{F}(s)=D_{k}(s)+M_{k}(s)+L_{k}(s),$$ 
 	where
 	
 	\begin{align}\label{DK}
 	D_{k}(s)=\int_{1}^{\infty}\bigg\{F(iy)-a_{0}\bigg\}y^{s-1}ds+i^{2k}\nu(T)\int_{1}^{\infty}\bigg\{F(iy)-a_{0} \bigg\}y^{2k-s-1}ds,
 	\end{align}
 	\begin{equation}\label{MK}
 	M_{k}(s)=a_{0}\bigg\{\frac{i^{2k}\nu(T)}{s-2k}-\frac{1}{s}\bigg\},
 	\end{equation}
 	and
 	
 	\begin{align}\label{LK}
 	L_{k} (s)& =  -\sum_{j=1}^{N} i^{2k+\alpha_{j}}\sum_{t=0}^{Mj} b_{jt}\sum_{l=0}^{t} {t\choose l} \left(\frac{i\pi}{2} \right)^{t-l} \frac{(-1)^{l+1}l!}{(2k+\alpha_{j}-s)^{l+1}}	,
 	\end{align}
 	for ~$\displaystyle Res>2k+\max_{j}|Re\alpha_{j}|.$
 \end{enumerate}
 \begin{enumerate}
 	\item[(B)] $\Phi_{F}(s)$~is bounded uniformly and absolutely in each~$\sigma$ in a  lacunary vertical strips of the form $$S(\sigma_{1},\sigma_{2},t_{0})=\lbrace s: \sigma_{1}\leq\sigma\leq\sigma_{2},|\Im s|=t\geq t_{0}\rbrace,$$ where~$\displaystyle\sigma_{1},\sigma_{2}\in\mathbb{R}$~and~$\displaystyle t_{0}>\max_{j}|\Im\alpha_{j}|.$
 \end{enumerate}
 Moreover, 
 \begin{enumerate}
 	\item[(C)]  $ \Phi_{F}(s)$ ~ satisfies the functional equation
 \end{enumerate}
 \begin{equation}\label{phi}
 \Phi_{F} \left( 2k-s\right) = e^{\pi ik}\nu(T)\Phi_{F}(s).
 \end{equation}	
 \end{thm}

{\em Proof.} Suppose ~$F$~is an entire automorphic integral of weight ~$2k,$~$k\in\mathbb{R},$~ multiplier system ~$\nu$ and associated ~LPPF~$q(z).$~For ~$Res>2k+\max_{j}|Re\alpha_{j}|,$ ~applying Mellin transform of $F,$ we have 
 \begin{eqnarray*}
 	\Phi_{F}(s)&=&\int_{0}^{\infty}\left( F(iy)-a_{0}\right)y^{s}\frac{dy}{y}\\
 	&=&\int_{1}^{\infty}\left( F(iy)-a_{0}\right)y^{s}\frac{dy}{y}+\int_{0}^{1}\left( F(iy)-a_{0}\right)y^{s}\frac{dy}{y}.
 \end{eqnarray*}
 Taking $y\mapsto\frac{1}{y}$ in the later integral and  applying the transformation law in (\ref{LPP2}),we get
 
 \begin{align*}
 \int_{0}^{1}\left( F(iy)-a_{0}\right)y^{s}\frac{dy}{y}&=\nu(T)i^{2k}\int_{1}^{\infty}\left( F(iy)-a_{0}\right)y^{2k-s-1}dy+\nu(T)i^{2k}a_{0}\int_{1}^{\infty}y^{2k-s-1}dy\\
 &+\nu(T)i^{2k}\int_{1}^{\infty}q(iy)y^{2k-s-1}dy-\frac{a_{0}}{s}\\
 &=\nu(T)i^{2k}\int_{1}^{\infty}\left( F(iy)-a_{0}\right)y^{2k-s-1}dy+a_{0}\bigg[\frac{\nu(T)i^{2k}}{s-2k}-\frac{1}{s}\bigg]\\
 &+\nu(T)i^{2k}\int_{1}^{\infty}q(iy)y^{2k-s-1}dy.
 \end{align*}
 
 Thus 
 \begin{align*}
 \Phi_{F}(s)=D_{k}(s)+M_{k}(s)+L_{k}(s),
 \end{align*}
 where
 \begin{align*}
 D_{k}(s)=\int_{1}^{\infty}\left( F(iy)-a_{0}\right)y^{s-1}dy+\nu(T)i^{2k}\int_{1}^{\infty}\left( F(iy)-a_{0}\right)y^{2k-s-1}dy,
 \end{align*}
 \begin{align*}
 M_{k}(s)=a_{0}\bigg[\frac{\nu(T)i^{2k}}{s-2k}-\frac{1}{s}\bigg],
 \end{align*}
 and 
 \begin{align*}
 L_{k}(s)=\nu(T)i^{2k}\int_{1}^{\infty}q(iy)y^{2k-s-1}dy.
 \end{align*}
 $D_{k}(s)$~is entire and converges uniformly (absolutely) on compact subset of ~$\mathbb{C},$~for large ~$Re (\ s)~~ and ~$  ~$D_{k}(s)$~satisfies  the functional equation $D_{k}(2k-s)=e^{\pi ik}\nu(T)D_{k}(s).$ Similarly $M_{k}(s)$~is meromorphic with simple poles at~$(s=2k,0),$~and satisfies the functional equation $M_{k}(2k-s)=e^{\pi ik}\nu(T)M_{k}(s).$ In order to show the meromerphic continuation of  $\Phi_{F}(s)$~we need the meromorphic  continuation of~ $L_{k}(s)$~to~$\mathbb{C}.$~To see this, note that
 \begin{align*}
 L_{k}(s)&=\nu(T)i^{2k}\int_{1}^{\infty}q(iy)y^{2k-s-1}dy\\
 &=\nu(T)i^{2k}\int_{1}^{\infty}\sum_{j=1}^{N}(iy)^{\alpha_{j}}\sum_{t=0}^{M_{j}}\beta_{jt}(\log iy)^{t}y^{2k-s-1}dy.
 \end{align*}
Using the fact that  $\log iy=\log y+i\arg(iy)=\log y+\frac{i\pi}{2}$ for~$y>0$, and applying the binomial theorem on $\left(\frac{i\pi}{2}+\log y \right)^{l}$, we can rewrite the above equation as

$$L_{k}(s)=\nu(T)\sum_{j=1}^{N}i^{2k+\alpha_{j}}\sum_{t=0}^{M_{j}}\beta_{jt}\sum_{l=0}^{t} {t\choose l}\left(\frac{i\pi}{2} \right)^{t-l}\int_{1}^{\infty}y^{2k+\alpha_{j}-s-1}(\log y)^{t}dy.$$
Integration by parts now yields 
$$ \displaystyle L_{k}(s)=\nu(T)\sum_{j=1}^{N}i^{2k+\alpha_{j}}\sum_{t=0}^{M_{j}}\beta_{jt}\sum_{l=0}^{t} {t\choose l}\left(\frac{i\pi}{2} \right)^{t-l} \frac{l!}{(s-(\alpha_{j}+2k))^{l+1}},$$ for $Re\ s>2k+\max_{j}Re\alpha_{j}.$ Thus $L_{k}(s)$ is meromorrhpic in $\mathbb{C}$ with poles of order $M_{j}+1$ at $ s=2k+\alpha_{j}.$ Thus 
$$\Phi_{F}(s)=D_{k}(s)+M_{k}(s)+L_{k}(s)$$ is meromorphic on $\mathbb{C}$ with finite simple poles at $s=2k,~ s=0 ~$  and poles of order $M_{j}+1$~at~ $s=2k+\alpha_{j}.$  This completes the prove of $(A).$ \\

To prove $\Phi_{F}(2k-s)=e^{\pi ik}\nu(T)\Phi(s)$ it is enough to show that $L_{k}(2k-s)=e^{\pi ik}\nu(T)L_{k}(s)$. For this we need  

\begin{lem}\label{EqnLppf}
	If $F$ and $q$ are as in Theorem \ref{ThM41}, then
	\begin{equation}
	q\left( \frac{-1}{z}\right)=-(-z)^{2k}\bar{\nu}(T)q(z), \forall z\in\mathcal{H}.
	\end{equation}
\end{lem}
 
{\em Proof.}  By taking~$z\rightarrow\frac{-1}{z}$~in (\ref{LPP2}), we have 
	$$ F(z)= \nu(T)(-z)^{-2k}\bigg[\nu(T)z^{2k}F(z)+\nu(T)z^{2k}q(z)+q\left(\frac{-1}{z}\right)\bigg].$$
	Thus,
	\begin{eqnarray*}
		F(z)-\nu^{2}(T)(-z)^{-2k}z^{2k}F(z)&=&\nu^{2}(-z)^{-2k}z^{2k}q(z)+\nu(T)(-z)^{-2k}q\left(\frac{-1}{z}\right)\\
		&=&i^{4k}\nu^{2}(T)q(z)+\nu(T)(-z)^{-2k}q\left(\frac{-1}{z}\right)\\
		&=&(i^{2}\nu(T))^{2}q(z)+\nu(T)(-z)^{-2k}q\left(\frac{-1}{z}\right)\\
		&=&q(z)+\nu(T)(-z)^{-2k}q\left(\frac{-1}{z}\right).
	\end{eqnarray*}
	Therefore, 
	$$ F(z)-\nu^{2}(T)(i^{2k})^{2}F(z)=q(z)+\nu(T)(-z)^{-2k}q\left(\frac{-1}{z}\right). $$
	From the consistency condition on $\nu$, we see that $\nu^{2}(T)(i^{2k})^{2}=1$ and hence  $$q(z)+\nu(T)(-z)^{-2k}q\left(\frac{-1}{z}\right)=0,~\forall z\in\mathcal{H}, $$ as desired.\\

		Now from~ $L_{k}(s)=\displaystyle i^{2k}\nu(T)\int_{1}^{\infty}y^{2k-s-1}q(iy)dy,$~~we have, 
	$$L_{k}(2k-s)=\displaystyle i^{2k-s}\nu(T)\int_{1}^{\infty}y^{s-1}q(iy)dy.$$   Applying Lemma \ref{EqnLppf} and  replace $y$ by $\frac{1}{y}$ in the integral we obtain $$\displaystyle i^{-2k}\bar{\nu}(T)L_{k}(2k-s)=-i^{2k}\nu(T)\int_{0}^{1}y^{2k-s-1}q(iy)dy.$$ \\
	Then by definition of $q$ we have  
	$$\displaystyle i^{-2k}\bar{\nu}(T)L_{k}(2k-s)=-i^{2k}\nu(T)\int_{0}^{1} \sum_{j=1}^{N}(iy)^{\alpha_{j}}\sum_{t=0}^{M_{j}}\beta_{jt}(\log iy)^{t}y^{2k-s-1}dy.$$\\
	
	Using  $\log iy=\frac{i\pi}{2}+\log y$ and applying the binomial theorem, the substitution  $\displaystyle y\mapsto \frac{1}{y}$, and integration by parts we have,
	
	 \begin{eqnarray*}\displaystyle i^{-2k}\bar{\nu}(T)L_{k}(2k-s)&=&-i^{2k}\nu(T) \sum_{j=1}^{N}(i)^{\alpha_{j}}\sum_{t=0}^{M_{j}}\beta_{jt} \sum_{l=0}^{t}\left(\frac{i\pi}{2} \right)^{t-l}\int_{1}^{\infty} y^{s-2k-\alpha_{j}-1}(\log y)^{l}(-1)^{l}dy\\
	 &=&\nu(T) \sum_{j=1}^{N}(i)^{\alpha_{j}+2k}\sum_{t=0}^{M_{j}}\beta_{jt} \sum_{l=0}^{t}{t\choose l} \left(\frac{i\pi}{2} \right)^{t-l} \frac{l!}{[s-(2k+\alpha_{j})]^{l+1}}.\end{eqnarray*} Hence  $$ L_{k}(s)=\nu(T) \sum_{j=1}^{N}(i)^{\alpha_{j}+2k}\sum_{t=0}^{M_{j}}\beta_{jt} \sum_{l=0}^{t}{t\choose l} \left(\frac{i\pi}{2} \right)^{t-l} \frac{l!}{[s-(2k+\alpha_{j})]^{l+1}}.$$
	   Therefore,  $\displaystyle L_{k}(2k-s)=e^{\pi ik}\nu(T)L_{k}(s)$.  
	This complete the prove of $(C).$ \\
	
	Finally to show the boundedness condition $(B)$, note that the rational functions $M_{k}$ and $L_{k}$ are bounded in lacunary vertical strips which do not contain poles of $M_{k}$ and $L_{k}$.  The boundedness of $D_{k}$ follows from the fact that $a_{m}=\mathcal{O}(m^{\beta})$ for $\sigma >1+\beta $ implies
	$$ \bigg|\int_{1}^{\infty} [F(iy)-a_{0}]y^{s-1}dy \bigg|=\mathcal{O}(\zeta(\sigma-\beta)\Gamma(\sigma)) . $$  This complete the proof of  $(B)$ and  that of the  theorem.\\

 Finally before we state and proof Theorem \ref{ThM42}, we state   Perron's Formula as the following Lemma. (See   \cite{CN} for details.)~   We shall use the convention of writing  $\displaystyle\int_{(b)}$ for   $\displaystyle\int_{b-i\infty}^{b+i\infty}$.
  
   \begin{lem}\label{LemPe}
  	Let    $\sigma_{0}$   be  the abscissa of absolute convergence for    $\displaystyle \varphi(s)=\sum_{m=1}^{\infty}a_{m}\lambda_{m}^{-s}$ and $\lbrace \lambda_{m} \rbrace$ be a  sequence of positive real numbers tending to $\infty$ as $m\rightarrow \infty.$ Then for    $k\geq 0,\sigma>0$   and   $\sigma>\sigma_{0}$, 
  	\begin{equation}
  	\frac{1}{\Gamma(k+1)}{\sum_{\lambda_{m}\leq x}}^{_\prime} a_{m}\left( x-\lambda_{m}\right)^{k}=\frac{1}{2\pi i}\int_{(\sigma)}\frac{\Gamma(s)\varphi(s)x^{s+k}}{\Gamma(s+k+1)}ds,
  	\end{equation}
  	where the prime $'$ on the summation sign indicates that if    $k=0$   and   $x=\lambda_{m}$   for some positive integer   $m$,    then we count only   $\frac{1}{2}a_{m}$.
  \end{lem}
 
 \section{Main results}
 
 In this section we state and prove the main results of the paper. 
   Here we establish two arithmetical identities associated to entire automorphic integrals involving log-polynomial period functions on the discrete Hecke group $G(\lambda),$ analogous to Chandrasekharan and Narasimhan and show that they are equivalent to the functional equation (\ref{phi}), $$\Phi_{F}(2k-s)=e^{\pi ik}\nu(T)\Phi_{F}(s).$$\\

\setcounter{equation}{0}

\begin{thm}\label{ThM42}
	Let $\Phi_{F}(s)$,$q(z)$ and $F(z)$ be as in Theorem~\ref{ThM41}. Then the functional equation 
\begin{equation} \label{THF}
\Phi_{F}(2k-s)=e^{\pi ik}\nu(T)\Phi_{F}(s)
\end{equation}
is equivalent to the arithmetical identity
\begin{align}\label{Efra2}
\frac{1}{\Gamma(\rho+1)}\sum_{0\leqslant m\leqslant x} a_{m}(x-m)^{\rho}&=e^{-\pi ik}\bar{\nu}(T)\left( \frac{2\pi}{\lambda}\right)^{-\rho}\sum_{m=1}^{\infty}a_{m} \left(\frac{x}{m} \right)^{\frac{2k+\rho}{2}}J_{2k+\rho}\left(\frac{4\pi \sqrt{mx}}{\lambda} \right) \nonumber \\
&+\left( \frac{2\pi}{\lambda}\right)^{2k}\frac{\nu(T) e^{\pi ik}a_{0} x^{2k+\rho}}{\Gamma(2k+\rho+1)}-\frac{a_{0} x^{\rho}}{\Gamma(\rho+1)} \nonumber\\
&+\left( \frac{2\pi}{\lambda}\right)^{2k+\alpha_{j}}\frac{x^{2k+\rho+\alpha_{j}}}{\Gamma(2k+\alpha_{j}+\rho+1)}\nu(T)e^{\pi ik}\sum_{j=1}^{N}i^{\alpha_{j}}\sum_{t=0}^{M_{j}}\beta_{jt}\left(\frac{i\pi}{2} \right)^{t},
\end{align}
for~$x>0,~\rho\in\mathbb{Z}_{\geq 0},~\alpha_{j},\beta_{jt}\in\mathbb{C},~\beta_{jM_{j}} \neq 0$~and~$N,M_{j}\in\mathbb{Z},$~provided~$\rho\geq 2\beta-2k,$~where~ $a_{m}=\mathcal{O}(m^{\beta}).$	
	
\end{thm}
 
 \subsection*{~~~~~~~~~~~~Proof of  $ (\ref{THF})$~ implies~$ (\ref{Efra2})$}.

  \subsection*{Proof}
  
  Applying Lemma \ref{LemPe}, for $\rho\geqslant 0, b\geqslant \beta$ and $a_{m}=\mathcal{O}(m^{\beta})$ we have  
  \begin{equation}\label{Perlppf}
  \frac{1}{\Gamma(\rho+1)}{\sum_{0\leqslant m\leqslant x}}^{\prime} a_{m}(x-m)^{\rho} =\frac{1}{2\pi i}\int_{(b)}\frac{\Gamma(s)\varphi(s)x^{s+\rho}}{\Gamma(s+\rho+1)}ds.
  \end{equation}
  Here the prime notation ' indicates that if $m=x, a_{m}$ is to  be multiplied by $\frac{1}{2}$ and $\displaystyle \varphi(s)= \sum_{m=1}^{\infty} \frac{a_{m}}{m^{s}} $ with $\displaystyle\sum_{m=1}^{\infty} \frac{|a_{m}|}{m^{\beta}}<\infty.$ \\
  To evaluate the integral on the right side of (\ref{Perlppf}), for $2k-b<\sigma<b,$ where $s=\sigma+it,$  consider rectangular region with vertices  $2k-b \pm iT$ and $b\pm iT,$ having positive orientation. Let $T$  be large positive  number  so that Stirling's formula can be used in approximating $\Gamma(s)$ as $|s|\rightarrow \infty.$ Using Stirling's formula and the Phragmen-Lindelof theorem we can show the integrals along the horizontal paths tend to $0$ as $T\rightarrow\infty.$ Hence with the substitution, $ \left(  \frac{2\pi}{\lambda}\right)^{s} \Phi_{F}(s)=\Gamma(s)\varphi(s),$  the right hand side of (\ref{Perlppf}) may be written as
  \begin{align*}
  \frac{1}{2\pi i} \int_{(b)} \frac{\Gamma(s)\varphi(s)}{\Gamma(s+\rho+1)}ds= & \frac{1}{2\pi i}\int_{(b)} \left(\frac{2\pi}{\lambda} \right)^{s}\frac{\Phi_{F}(s) x^{s+\rho}}{\Gamma(s+\rho+1)} ds\\
  & =\frac{1}{2\pi i}\int_{(2k-b)} \left(\frac{2\pi}{\lambda} \right)^{s}\frac{\Phi_{F}(s) x^{s+\rho}}{\Gamma(s+\rho+1)} ds\\
  & +\sum_{s\in Polset of \Phi_{F}(s)}Res\bigg\lbrace \left(\frac{2\pi}{\lambda} \right)^{s}\frac{\Phi_{F}(s) x^{s+\rho}}{\Gamma(s+\rho+1)}\bigg\rbrace .
  \end{align*}
  Put $$\omega_{1}(x,\lambda)=\frac{1}{2\pi i}\int_{(2k-b)} \left(\frac{2\pi}{\lambda} \right)^{s}\frac{\Phi_{F}(s) x^{s+\rho}}{\Gamma(s+\rho+1)} ds $$ and
  $$\omega_{2}(x,\lambda)=\sum_{s\in Polset of \Phi_{F}(s)}Res\bigg\lbrace \left(\frac{2\pi}{\lambda} \right)^{s}\frac{\Phi_{F}(s) x^{s+\rho}}{\Gamma(s+\rho+1)}\bigg\rbrace .$$ To evaluate $ \omega_{1}(x,\lambda)$ we begin  by using the functional equation 
  $$ \Phi_{F}(2k-s)=e^{\pi ik}\nu(T)\Phi(s).$$ 
  $\omega_{1}(x,\lambda)$ can be written as 
  $$\omega_{1}(x,\lambda)=\frac{1}{2\pi i}\int_{(2k-b)} \left(\frac{2\pi}{\lambda} \right)^{s}\frac{e^{-\pi ik}\bar{\nu}(T) \Phi_{F}(2k-s) x^{s+\rho}}{\Gamma(s+\rho+1)} ds. $$
  By substituting $$\Phi_{F}(2k-s)= \left(\frac{2\pi}{\lambda} \right)^{s}\Gamma(2k-s)\sum_{m=1}^{\infty}\frac{a_{m}}{m^{2k-s}}$$ and replacing $s$ by~$2k-s,$  to evaluate~ $\omega_{1}(x,\lambda)$~ we can use the integral in  (\cite{CN},~~page $5$ )  given  by the formula 
  
  \begin{equation*}
  \frac{1}{2\pi i}\int_{C-i\infty}^{C+i\infty}\frac{2^{s-\nu-1 \Gamma(\frac{s}{2})}}{\Gamma(\nu-\frac{1}{2}s+1)} x^{-s} ds=\frac{J_{\nu}(x)}{x^{\nu}},
  \end{equation*}
  provided that   $0<C\leq \nu+1,~~   \nu>0.$
  
   With simple algebraic manipulation we obtain 
  \begin{equation}\label{eq1lppf}
  \omega_{1}(x,\lambda)=\left(\frac{2\pi}{\lambda} \right)^{-\rho}e^{-\pi ik}\bar{\nu}(T)\sum_{m=1}^{\infty}a_{m} \left( \frac{x}{m}\right)^{\frac{2k+\rho}{2}}J_{2k+\rho}\left(\frac{4\pi\sqrt{mx}}{\lambda} \right),
  \end{equation}
  provided that~$\rho \geq 2\beta-2k.$ \\
  To complete  evaluation of the integral in (\ref{Perlppf}) we calculate the residues  in $\omega_{2}(x,\lambda).$  Note that $\Phi_{F}$~ has poles at~$s=0,2k$~ resulting from $M_{k}$ and at $s=2k+\alpha_{j}$ with order $M_{j}+1$ from $L_{k}$.  Thus  calculating the residues of the functions we obtain,
  \begin{align}\label{eq2lppf}
  \omega_{2}(x,\lambda)=&\left(\frac{2\pi}{\lambda} \right)^{2k} e^{\pi ik}\frac{x^{2k+\rho}a_{0}\nu(T)}{\Gamma(2k+\rho+1)}-\frac{a_{0}x^{\rho}}{\Gamma(\rho+1)} \nonumber \\ 
  & +\left(\frac{2\pi}{\lambda} \right)^{2k+\alpha_{j}}\frac{x^{2k+\rho+\alpha_{j}}}{\Gamma(2k+\rho+1+\alpha_{j})}\nu(T)e^{\pi ik}\sum_{j=1}^{N}i^{\alpha_{j}}\sum_{t=0}^{M_{j}}\beta_{jt}\left(\frac{i\pi}{2} \right)^{t},
  \end{align}
  for~$\rho\in\mathbb{Z}_{\geq 0},$~$\alpha_{j},\beta_{jt}\in\mathbb{C}$~and~$\beta_{jM_{j}}\neq  0.$\\
  The integral in  right of (\ref{Perlppf}) is equal to the combination of the expressions in (\ref{eq1lppf}) and (\ref{eq2lppf}) respectively, thus it follows the identity (\ref{Efra2}).

  \subsection*{~Proof of converse of the theorem:}
  
 Suppose $F$ be entire automorphic integral function. $F(z)$ has Fourier series expansion of the form 
 \begin{equation}\label{fouri}
 F(z)=\sum_{m=0}^{\infty}a_{m}e^{2\pi im/\lambda}.
 \end{equation}
 For $z\in\mathcal{H},$ where $a_{m}=\mathcal{O}(m^{\gamma})$ as $m\rightarrow\infty,$ for some $\gamma>0.$  Put $$\Phi_{F}(s)=\left( \frac{2\pi}{\lambda}\right)^{-s}\Gamma(s)\sum_{m=1}^{\infty}a_{m}m^{-s}.$$  
 For $z\in\mathcal{H},$~ $F$ satisfies the relation
 \begin{equation}\label{trans}
 \bar{\nu}(T)z^{-2k}F\left(\frac{-1}{z} \right)=F(z)+q(z)\forall z\in\mathcal{H},
 \end{equation}
 where $$q(z)=\sum_{j=1}^{N}z^{\alpha_{j}}\sum_{t=0}^{M_{j}}\beta_{jt}(\log z)^{t}, ~~\alpha_{j}, \beta_{jt}\in\mathbb{C}.$$
 Now using  (\ref{fouri}) and (\ref{trans}) we have 
 $$ z^{-2k}\bar{\nu}(T)\sum_{m=0}^{\infty}a_{m}e^{-2\pi im/\lambda z}=\sum_{m=0}^{\infty}a_{m}e^{2\pi imz/\lambda}+q(z).$$ Letting $z=\frac{iy\lambda}{2\pi}, y>c, c\in\mathbb{R^{+}},$~   we have 
 \begin{equation}\label{lppfy}
 \left(\frac{2\pi}{iy\lambda} \right)^{2k}\bar{\nu}(T)\sum_{m=0}^{\infty}a_{m}e^{\frac{-4\pi^{2}m}{y\lambda^{2}}}=\sum_{m=0}^{\infty}a_{m}e^{-my}+q\left(\frac{iy\lambda}{2\pi} \right).
 \end{equation}
 In proving the converse, it is enough to show that the identity in (\ref{Efra2}) implies (\ref{lppfy}). Since $q(z)$ is entire in $\mathcal{H}$ and the series expansion for $F$ in (\ref{fouri}) is uniformly convergent in compact subsets of $\mathcal{H},$ the automorphic relation in (\ref{trans}) will follow by analytic for $\forall z\in\mathcal{H}.$   Following the method used in Chandraseharan and Narasimhan \cite{CN},  we multiply the identity in (\ref{Efra2}) through out by $y^{\rho+1}e^{-xy}$ with $y>0,$ integrate relative to $x$ from $0$ to $\infty.$
 To this end consider four integrals (separate) corresponding to the terms occurring in (\ref{Efra2}) as 
 \begin{align*}
 \varpi_{1}(y,\lambda)&=\int_{0}^{\infty}\Bigg\lbrace\frac{1}{\Gamma(\rho+1)}\sum_{0\leq m\leq x} a_{m}(x-m)^{\rho}\Bigg\rbrace y^{\rho+1}e^{-xy}dx.
 \end{align*}
 \begin{align*}
 \varpi_{2}(y,\lambda)&= \int_{0}^{\infty}\Bigg\lbrace e^{-\pi ik}\bar{\nu}(T) \left(\frac{2\pi}{\lambda} \right)^{-\rho}\sum_{m=1}^{\infty}a_{m}\left(\frac{x}{m} \right)^{\frac{2k+\rho}{2}}J_{2k+\rho}\left(\frac{4\pi\sqrt{mx}}{\lambda} \right) \Bigg\rbrace y^{\rho+1}e^{-xy}dx.
 \end{align*}
 \begin{align*}
 \varpi_{3}(y,\lambda)&= \int_{0}^{\infty} \Bigg\lbrace e^{\pi ik}\left(\frac{2\pi}{\lambda} \right)^{2k} \frac{a_{0}\nu(T)x^{2k+\rho}}{\Gamma(2k+\rho+1)}-\frac{a_{0}x^{\rho}}{\Gamma(\rho+1)}\bigg\rbrace y^{\rho+1}e^{-xy}dx.
 \end{align*}
 \begin{align*}
 \varpi_{4}(y,\lambda)&=\int_{0}^{\infty} \Bigg\lbrace\left(\frac{2\pi}{\lambda} \right)^{2k+\alpha_{j}} \frac{x^{2k+\rho+\alpha_{j}}}{\Gamma(2k+\rho+\alpha_{j}+1)}\nu(T) e^{\pi ik}\sum_{j=1}^{N}i^{\alpha_{j}}\sum_{t=0}^{M_{j}}\beta_{jt}\left( \frac{i\pi}{2}\right)^{t}\Bigg\rbrace  y^{\rho+1}e^{-xy}dx.
 \end{align*}
 The evaluation of $\varpi_{1}$ and $\varpi_{2}$  follows as Chandrasekharan and Narasimhan found in \cite{CN}.  Interchange of summation and integration being permitted for $\rho \geq 2\beta-2k-\frac{1}{2}.$ Then
 \begin{align*}
 \varpi_{1}(y,\lambda)=\frac{y^{\rho}+1}{\Gamma(\rho+1)}\int_{0}^{\infty}x^{\rho}e^{-xy}dx+\frac{y^{\rho}+1}{\Gamma(\rho+1)}\sum_{1\leq m\leq x}a_{m}\int_{m}^{\infty}(x-\rho)^{\rho}e^{-xy}dx.
 \end{align*}
 Integration by substitution and using the standard integral representation for $\Gamma(s)$ we get
 \begin{equation}\label{eqq1}
 \varpi_{1}(y,\lambda)=\sum_{m=0}^{\infty}a_{m}e^{-my}.
 \end{equation}
 To compute $\varpi_{2}$ we apply the  integral formula $\int_{0}^{\infty}e^{-xy}J_{\nu}(a\sqrt{x})x^{\frac{\nu}{2}}dx=\frac{2a^{\nu}}{(2y)^{\nu+1}}e^{\frac{-a^2}{4y}}$, provided that $y>0, Re(\nu)>-1,a>0.$ And  with the substitution ~$\frac{4\pi\sqrt{m}}{\lambda}=a, \nu=\rho+2k$~ after simplifying, we obtain
 \begin{equation}\label{eqq2}
 \varpi_{2}(y,\lambda)=\left(\frac{2\pi}{\lambda y} \right)^{2k}e^{-\pi ik}\bar{\nu}(T)\sum_{m=1}^{\infty}a_{m}e^{\frac{-4\pi^{2} m}{y\lambda^{2}}}.
 \end{equation}
 The integral in~ $ \varpi_{3}(y,\lambda)$ is evaluated  using the integral representation for $\Gamma(s).$  Therefore. after simplifying, we see that
 \begin{equation}\label{eqq3}
 \varpi_{3}(y,\lambda)=\left(\frac{2\pi}{\lambda y} \right)^{2k} a_{0}e^{\pi ik}\nu(T)-a_{0}.
 \end{equation}
 To evaluate $\varpi_{4}(y,\lambda)$ since each term is integrable the inter change of the integration and the finite double sum is valid . Then integration by substitution  and using the integral representation of $\Gamma(z)$  we obtain
 \begin{align}\label{eqq4}
 \varpi_{4}(y,\lambda)&=\left(\frac{2\pi i}{\lambda y} \right)^{2k}\nu(T)\sum_{j=1}^{N}\left( \frac{2\pi i}{\lambda y} \right)^{\alpha_{j}}\sum_{t=0}^{M_{j}}\beta_{jt}\left( \frac{i\pi}{2}\right)^{t}\nonumber \\ 
 & =\left(\frac{-1}{i\lambda y/2\pi} \right)^{2k}\nu(T)\sum_{j=1}^{N}\left( \frac{-1}{\lambda yi/2\pi} \right)^{\alpha_{j}}\sum_{t=0}^{M_{j}}\beta_{jt}\left( \frac{i\pi}{2}\right)^{t}.
 \end{align}
 For $y>0,$ we see that $\log\left(\frac{-1}{iy} \right)+\log y=\frac{i\pi}{2}.$ Then for $l\in\mathbb{N}\bigcup\lbrace 0\rbrace$ and  applying the binomial theorem we have 
 $$\left( \log\left(\frac{-1}{iy} \right)+\log y \right)^{t}=\sum_{l=0}^{t} {t\choose l} \left(\log\left(\frac{-1}{iy} \right)\right)^{t-l} (\log y)^{l}.$$
 If $l=0$ we have $ \left(\log \left(\frac{-1}{iy} \right)\right)^{t}=\left( \log\left(\frac{-1}{iy} \right)+\log y \right)^{t}. $ Now replacing~$\frac{-1}{iy} $~by~$\frac{-1}{i\lambda y/2\pi} $~and substituting for~ $\left( \frac{i\pi}{2}\right)^{t} $~in~ (\ref{eqq4}) yields 
 \begin{equation}
 \varpi_{4}(y,\lambda) 
 =\left(\frac{-1}{i\lambda y/2\pi} \right)^{2k}\nu(T)\sum_{j=1}^{N}\left( \frac{-1}{\lambda yi/2\pi} \right)^{\alpha_{j}}\sum_{t=0}^{M_{j}}\beta_{jt} \left(\log\left(\frac{-1}{iy\lambda/2\pi} \right)\right)^{t}.
 \end{equation}
 Observe that $\varpi_{4}(y,\lambda)=\left(\frac{-1}{i\lambda y/2\pi} \right)^{2k}\nu(T)q\left(\frac{-1}{iy\lambda/2\pi} \right),$ where $q(z)$ is the log-polynomial period function in (\ref{trans}). Then applying Lemma\ref{EqnLppf} we have $\varpi_{4}(x,\lambda) =\left(\frac{-1}{i\lambda y/2\pi} \right)^{2k}\nu(T)q\left(\frac{-1}{iy\lambda/2\pi} \right)=-q\left(\frac{iy\lambda}{2\pi} \right).$
 Combining the results in (\ref{eqq1}), (\ref{eqq2}), (\ref{eqq3}) and (\ref{eqq4}) we showed that the identity in (\ref{Efra2}) implies
 \begin{align}\label{Finee}
 \sum_{m=0}^{\infty}a_{m}e^{-my}=\left(\frac{2\pi}{\lambda y} \right)^{2k}e^{-\pi ik}\bar{\nu}(T)\sum_{m=1}^{\infty}a_{m}e^{\frac{-4\pi^{2} m}{y\lambda^{2}}}
 +\left(\frac{2\pi}{\lambda y} \right)^{2k} a_{0}e^{\pi ik}\nu(T) \nonumber \\ 
 +\left(\frac{-1}{i\lambda y/2\pi} \right)^{2k}\nu(T)\sum_{j=1}^{N}\left( \frac{-1}{\lambda yi/2\pi} \right)^{\alpha_{j}}\sum_{t=0}^{M_{j}}\beta_{jt}\left(\log\left( \frac{-1}{iy\lambda/2\pi}\right) \right)^{t}.
 \end{align}
 Observe that $$\left(\frac{-1}{i\lambda y/2\pi} \right)^{2k}\nu(T)\sum_{j=1}^{N}\left( \frac{-1}{\lambda yi/2\pi} \right)^{\alpha_{j}}\sum_{t=0}^{M_{j}}\beta_{jt}\left(\log\left( \frac{-1}{iy\lambda/2\pi}\right) \right)^{t}=-q\left(\frac{iy\lambda}{2\pi} \right).$$ 
 Since  $F(z)$ is holomorphic and has a Fourier series expansion $\displaystyle F(z)=\sum_{m=0}^{\infty} a_{m}e^{2\pi imz/\lambda}$ in $\mathcal{H}.$ Then  with $z=\frac{iy\lambda}{2\pi}$ (\ref{Finee}) written as (\ref{lppfy}). By the identity theorem ,then the automorphic relation (\ref{LPP2}) follows $\forall z\in\mathcal{H},$ and this completes  proof of the theorem.
  

\maketitle
\def\theequation{\thesection.\arabic{equation}}
 
\setcounter{equation}{0}

 \begin{thm}\label{ThM43}
 	Let $ \Phi_{F}(s), F(z)$ and $q(z)$ be as in  Theorem\ref{ThM41} . Then the functional equation 
 \begin{equation}\label{Func2}
 \Phi_{F}(2k-s)=e^{\pi ik}\nu(T)\Phi_{F}(s)
 \end{equation} 
 is equivalent to the arithmetical identity 
 \begin{align}\label{Efrata2}
 \left( \frac{-1}{y} \frac{d}{dy}\right)^{\rho}\left(\frac{1}{y}\sum_{m=1}^{\infty}a_{m}e^{-y\sqrt{m}} \right)& = \frac{2^{\rho}}{\sqrt{\pi}} \Gamma\left(2k+\rho+\frac{1}{2}\right)\left( \frac{8\pi}{\lambda}\right)^{2k} \bar{\nu}(T) e^{-\pi ik}\nonumber\\
 &\times\sum_{m=0}^{\infty}\frac{a_{m}}{\left(y^{2}+4\left(\frac{2\pi}{\lambda}\right)^{2}m \right)^{2k+\rho+\frac{1}{2}}}\nonumber \\
 & +\frac{2^{\rho}}{\sqrt{\pi}}\frac{\Gamma(2k+\rho+\alpha_{j}+\frac{1}{2})}{y^{2\rho+1}}\left( \frac{8\pi}{\lambda y^{2}}\right)^{2k} \nu(T) e^{\pi ik}\nonumber\\
 &\times\sum_{j=1}^{N}\left(\frac{8\pi i}{\lambda y^{2}} \right)^{\alpha_{j}}\sum_{t=0}^{M_{j}}\beta_{jt}\left(\frac{i\pi}{2} \right)^{t} \nonumber \\
 & -\frac{2^{\rho+1}}{\sqrt{\pi}}\frac{a_{0}}{y^{2\rho+1}} \Gamma\left(\rho+\frac{1}{2}\right),
 \end{align}
 provided that~Re$y\in\mathbb{R^{+}},~\alpha_{j},\beta_{jt}\in\mathbb{C},~\rho\in\mathbb{Z}_{\geq 0},$~and~$\rho >\beta-2k+\frac{1}{2},$\\
 ~where~~$\sum_{m=1}^{\infty} \frac{|a_{m}|}{m^{\beta}} <\infty$.	
 \end{thm}

Since, by Theorem \ref{ThM42} the functional equation (\ref{Func2}) is equivalent to (\ref{Efra2}), it would be enough to show  that the functional equation (\ref{Func2}) implies (\ref{Efrata2}) and that (\ref{Efrata2}) in turn implies the identity (\ref{Efra2}).
 
 \subsection*{Proof of implication:}
 
 	First we begin by showing that  (\ref{Func2}) implies (\ref{Efrata2}). Let $ \rho$ be an integer, $\rho \geq 0$ and  $\displaystyle \varphi(s)=\sum_{m=1}^{\infty}a_{m} m^{-s}$ with   $a_{m}=\mathcal{O}(m^{\beta}).$ Then for $\gamma>0, \gamma\geq \beta$ and applying the version of Perron's in Lemma  \ref{LemPe},  we have
 \begin{equation}\label{PerLppf}
 \frac{1}{\Gamma(\rho+1)}{\sum_{0\leq m \leq x}}^\prime a_{m}(x-m)^\rho = \frac{1}{2\pi i}\int_{(\gamma)} \frac{\varphi(s)\Gamma(s)x^{s+\rho}}{\Gamma(s+\rho+1)}ds.
 \end{equation} 
 As Chandrasekharan and Narasimhan in \cite{CN} we multiply (\ref{PerLppf})  by $e^{-y\sqrt{x}}x^{\frac{1}{2}}$ and integrating with respect to the variable $x$ on $[~0,\infty).$  For $\gamma >2k,$ Re$s=\gamma$ is the vertical path of integration. 
 Thus, we need to prove the identity (\ref{Efrata2}) follows from the equation
 \begin{multline}\label{MulT}
 \int^\infty_0e^{-y\sqrt{x}}
 x^{-\frac{1}{2}}\Bigg\{\frac{1}{\Gamma(\rho+1)}{\sum_{0\leq m\leq
 		x}}^\prime a_{m}{(x-m)}^\rho\Bigg\}dx\\
 =\int^\infty_0e^{-y\sqrt{x}}x^{-\frac{1}{2}}\Bigg\{\frac{1}{2\pi i}
 \int_{(\gamma)}
 \frac{\varphi(s)\Gamma(s)x^{s+\rho}}{\Gamma(s+\rho+1)}ds\Bigg\}dx,
 \end{multline}
 where $\rho\in\mathbb{Z}_{\geq 0},
 \rho\geq\gamma-2k+\frac{1}{2},$ and $y\in\mathbb{R}^+.$ Chandrasekharan and Narasimhan \cite{CN} for $\lambda_{n}\in\mathbb{R^{+}}$~$\lambda_{n}\rightarrow \infty$~as~$n\rightarrow\infty$ showed the identity
 \begin{equation}\label{chan}
 \sum^\infty_{n=1} a_{n} \int_{\lambda_n}^{\infty}\frac{(x-n)^\rho}{\Gamma(\rho+1)}e^{-y\sqrt{x}}x^{-\frac{1}{2}}dx
 =2(-2)^\rho\left(\frac{1}{y}\frac{d}{dy}\right)^\rho\left[\frac{1}{y}\sum^\infty_{n=1}a_ne^{-y\sqrt{\lambda_n}}\right].
 \end{equation}
 Then by the identity in (\ref{chan}) with $\lambda_{n}=m,$ we have 
 $$\sum^\infty_{m=1} a_{m} \int_{m}^{\infty}\frac{(x-m)^\rho}{\Gamma(\rho+1)}e^{-y\sqrt{x}}x^{-\frac{1}{2}}dx
 =2(2)^\rho\left(\frac{-1}{y}\frac{d}{dy}\right)^\rho\left[\frac{1}{y}\sum^\infty_{m=1}a_me^{-y\sqrt{m}}\right]. $$   Now on the left side of (\ref{MulT}) the interchange of integration and summation is justified for $\rho\geq 0.$ Then 
 \begin{align*}
 \int^\infty_0e^{-y\sqrt{x}}
 x^{-\frac{1}{2}}\Bigg\{\frac{1}{\Gamma(\rho+1)}{\sum_{0\leq m\leq
 		x}}^\prime a_{m}{(x-m)}^\rho\Bigg\}dx&=\frac{a_{0}}{\Gamma(\rho+1)} \int_{0}^{\infty}e^{-y\sqrt{x}}x^{\rho-\frac{1}{2}}dx\\
 &+ \frac{1}{\Gamma(\rho+1)}\times\\ 
 &\sum_{1\leq m\leq x} a_{m}
 \int^\infty_m e^{-y\sqrt{x}}x^{-\frac{1}{2}}(x-m)^{\rho}dx .
 \end{align*}
 Thus applying (\ref{chan}) we have 
 
 \begin{align}
 \int^\infty_0e^{-y\sqrt{x}}
 x^{-\frac{1}{2}}\Bigg\{\frac{1}{\Gamma(\rho+1)}{\sum_{0\leq m\leq
 		x}}^\prime a_{m}{(x-m)}^\rho\Bigg\}dx&=\frac{2^{2\rho+1}}{y^{2\rho+1}}a_{0}\frac{\Gamma(\rho+\frac{1}{2})}{\sqrt{\pi}}\nonumber\\
 &+2^{\rho+1}\left(\frac{-1}{y}\frac{d}{dp} \right)^{\rho}\bigg[ \frac{1}{y}\sum_{m=1}^{\infty} a_{m}e^{y\sqrt{m}}\bigg].
 \end{align}
 
 Now assuming $\Phi_{F}(s)= \left( \frac{2\pi}{\lambda}\right)^{-s} \Gamma(s)\varphi(s),$  where Re$s>\beta$ the right-hand side of (\ref{MulT}) may be written as 
 $$ W(y,\lambda)=\int^\infty_0e^{-y\sqrt{x}}x^{-\frac{1}{2}}\Bigg\{\frac{1}{2\pi i}
 \int_{(\gamma)}
 \left( \frac{2\pi}{\lambda}\right)^{s}\frac{ \Phi_{F}(s) x^{s+\rho}}{\Gamma(s+\rho+1)}ds\Bigg\}dx ,$$ for $\gamma \geq \beta.$ We may interchange the order of integration for $\rho \geq 0.$ 
 Then we have  $$W(y,\lambda)=\frac{1}{2\pi i}
 \int_{(\gamma)}
 \Bigg\lbrace \left( \frac{2\pi}{\lambda}\right)^{s}\frac{ \Phi_{F}(s)}{\Gamma(s+\rho+1)} \int^\infty_0e^{-y\sqrt{x}}x^{s+\rho-\frac{1}{2}} dx\Bigg\rbrace ds. $$ Using integration by substitution and applying the Legendre duplication formula  
 \begin{equation}\label{PR3}
 \sqrt{\pi} \Gamma(2z)=2^{2z-1}\Gamma(z)\Gamma\left(z+\frac{1}{2}\right),
 \end{equation}
  for~ ~$\Gamma(2s+2\rho+1),$ we obtain 
 $$ W(y,\lambda)=\frac{1}{2\pi i}
 \int_{(\gamma)}
 \left( \frac{2\pi}{\lambda}\right)^{s} \Phi_{F}(s)\frac{2^{2\rho+2s+1}}{y^{2\rho+2s+1}} \frac{\Gamma(s+\rho+\frac{1}{2})}{\sqrt{\pi}}ds .$$  To evaluate this integral we shall consider a rectangular region with vertices $2k-\gamma\pm iT$ and $\gamma \pm iT,$  $s=\sigma+it$ on its interior.  Using Cauchy's  residue theorem, the line of integration oriented positively to  Re$s=2k-\gamma.$ If $\rho \geq \gamma-2k+\frac{1}{2}$ the poles of $\Gamma\left( s+\rho+\frac{1}{2}\right)$ lie to the left of the vertical path. Then the poles of the integrand in the interval $[2k-\gamma, \gamma]$ arise from the function $\Phi_{F}(s).$ Applying Stirling's formula and the Phragmen-Lindelof theorem the integrals along both horizontal paths tend to $0$ as $T\rightarrow \infty,$   and thus  $ W(y,\lambda)$ may be written as 
 $$  W(y,\lambda) =w_{1} (y,\lambda)+w_{2}(y,\lambda),$$ where 
 \begin{align*}
 w_{1}(y,\lambda)=\frac{1}{2\pi i}\int_{(2k-\gamma)}\left(\frac{2\pi}{\lambda} \right)^{s}\Phi_{F}(s)\frac{\Gamma(s+\rho+\frac{1}{2})}{\sqrt{\pi}}\frac{2^{2s+2\rho+1}}{y^{2s+2\rho+1}}ds,
 \end{align*}
 and
 \begin{align*}
 w_{2}(y,\lambda)=\sum_{s\in polset\Phi(s)}\text{Res} \bigg\lbrace \left(\frac{2\pi}{\lambda} \right)^{s}\Phi(s)\frac{\Gamma(s+\rho+\frac{1}{2})}{\sqrt{\pi}}\frac{2^{2s+2\rho+1}}{y^{2s+2\rho+1}}\bigg\rbrace.
 \end{align*}
 To evaluate  $w_{1}(y,\lambda)$  recall the functional equation $$ \Phi_{F}(2k-s)=e^{\pi ik}\nu(T)\Phi_{F}(s)$$ and substitute in to the integral of $w_{1}(s,\lambda).$
 Then we have $$ w_{1}(y,\lambda)=\bar{\nu}(T)e^{-\pi ik} \frac{2^{2\rho+1}}{y^{2\rho+1}}\Bigg\lbrace  \frac{1}{2\pi i} \int_{(2k-\gamma)}\left(\frac{2\pi}{\lambda}\right)^{s} \Phi_{F}(2k-s) \frac{2^{2s}\Gamma(s+\rho+\frac{1}{2})}{y^{2s}\sqrt{\pi}}  ds\Bigg\rbrace.$$ Now replace~$\theta$~by~$2k-s$~and then setting~ $\Phi_{F}(\theta)= \left(\frac{2\pi}{\lambda} \right)^{-\theta}\Gamma(\theta)\sum_{m=1}^{\infty} a_{m}m^{-\theta}.$~By organizing the expressions~ $w_{1}(y,\lambda)$~may be written as 
 $$w_{1}(y,\lambda)= \bar{\nu}(T)e^{-\pi ik} \frac{2^{2\rho+1}}{y^{2\rho+1}}\sum_{m=1}^{\infty}a_{m}\Bigg\lbrace \frac{1}{2\pi i} \int_{(\gamma)}\Gamma(\theta)\Gamma(2k-\theta+\rho+\frac{1}{2})\left(\frac{16m\pi^{2}}{\lambda^{2}y^{2}} \right)^{-\theta} d\theta \Bigg\rbrace .$$ 
 The interchange of the integration with the summation is justified for~$\rho>\beta-2k+\frac{1}{2}.$\\ Replace~$\theta$~by~$-\theta$~and using the integral formula 
 \begin{equation}\label{PR4}
 \frac{1}{2\pi i}\int_{\gamma-i\infty}^{\gamma+i\infty}\Gamma(-s)\Gamma(\beta+s)t^{s}ds=\Gamma(\beta)(1+t)^{-\beta},
 \end{equation}
 where    $0>\gamma>Re(1-\beta)$   and $|arg   t|<\pi .$
  After simplifying, we obtain 
 \begin{equation}\label{w1}
 w_{1}(y,\lambda)=\frac{\bar{\nu}(T)e^{-\pi ik}}{\sqrt{\pi}} 2^{4k+2\rho+1}\Gamma\left(2k+\rho+\frac{1}{2} \right)\left( \frac{2\pi}{\lambda}\right)^{2k} \sum_{m=1}^{\infty} \frac{a_{m}}{\left( y^{2}+(\frac{4\pi}{\lambda})^{2}m\right)^{2k+\rho+\frac{1}{2}}}.
 \end{equation} 
 The series on~ (\ref{w1})~converges absolutely for~Re$y\in\mathbb{R^{+}},$~provided~$\rho >\beta-2k+\frac{1}{2}.$\\
 To proceed evaluation of  $ w_{2}(y,\lambda)$  recall that in Theorem \ref{ThM41} we have that $\Phi_{F}(s)=D_{k}(s)+M_{k}(s)+L_{k}(s),$ where $D_{k}(s),M_{k}(s)$ and $L_{k}(s)$ are functions in (\ref{DK}), (\ref{MK}), and (\ref{LK}) respectively. For $\rho\geq \gamma-2k+\frac{1}{2} $  the function $\left(\frac{2\pi}{\lambda} \right)^{s}\frac{\Gamma(s+\rho+\frac{1}{2})}{\sqrt{\pi}} \frac{2^{2s+2\rho+1}}{y^{2s+2\rho+1}}$ is analytic in the region bounded by Re$s=2k-\gamma$ and Re$s=\gamma,$  thus the poles of the intgrand in the interval $[2k-\gamma, \gamma]$ arise from  $\Phi_{F}(s).$  We also note that $D_{k}(s)$ is entire, $M_{k}(s)$ has simple poles at $(s=2k, 0)$ and $L_{k}(s)$ has a pole at $s=2k+\alpha_{j}$ of order $M_{j}+1,$ where $\alpha_{j}\in\mathbb{C}, M_{j}\in\mathbb{Z}_{\geq 0} ~.$  Therefore. by computing the residues of $M_{k}$ and $ L_{k}$ after simplifying, we obtain
 \begin{align}\label{W2}
 w_{2}(y,\lambda)=\frac{2^{2\rho+1}}{\sqrt{\pi}y^{\rho+1}}a_{0}\nu(T)e^{\pi ik}\left(\frac{8\pi}{\lambda y^{2}} \right)^{2k}\Gamma\left( 2k+\rho+\frac{1}{2}\right)-\frac{2^{2\rho+1}}{\sqrt{\pi}y^{\rho+1}}a_{0}\Gamma\left(\rho+\frac{1}{2}\right)\nonumber \\ 
 +\frac{2^{2\rho+1}}{\sqrt{\pi}y^{\rho+1}} \nu(T)e^{\pi ik} \left(\frac{8\pi}{\lambda y^{2}} \right)^{2k}\Gamma\left( 2k+\rho+\alpha_{j}+\frac{1}{2}\right)\sum_{j=1}^{N}i^{\alpha_{j}}\sum_{t=0}^{M_{j}}\beta_{jt}\left( \frac{i\pi}{2}\right)^{t}.
 \end{align} 
 There for by rewriting $W(y,\lambda)$ as combination of (\ref{w1}) and (\ref{W2}) and substituting the respective expressions in to the right-hand side of (\ref{MulT}), with simple rearrangement, we obtain the identity  (\ref{Efrata2}), provided that $\rho\in\mathbb{Z}_{\geq 0},~\rho \geq \beta-2k+\frac{1}{2}$ and $y\in\mathbb{R^{+}},$~where ~$\sum_{m=1}^{\infty} \frac{|a_{m}|}{m^{\beta}} <\infty.$

 \subsection*{Proof of the converse:}
 
 	To prove the converse it suffices to show that (\ref{Efrata2}) implies (\ref{Efra2}). As technique of Chandrasekharan and Narasimhan in \cite{CN}  multiply  (\ref{Efrata2}) by $e^{y\sqrt{x}},$  with Re$y>0,$ and $x>0$ and integrate the expression along a vertical path Re$s=\theta,$ where $\theta >0.$ The left hand side of the transformed equation can be evaluated using the formula
 \begin{align}
 \sum_{m=1}^{\infty} a_{m}\frac{1}{2\pi i}\int_{(\theta)}e^{y\sqrt{x}}\left(\frac{-1}{y}\frac{d}{dy} \right)^{\rho}\bigg[\frac{1}{y}e^{-y\sqrt{m}}\bigg]dy=\frac{1}{\Gamma(\rho+1)}{\sum_{m\leq x}}^{\prime}(x-m)^{\rho}2^{-\rho},
 \end{align}
 while the right hand side of (\ref{Efrata2}) we calculate the integral of each terms one by one . So put
 \begin{align}
 h_{1}(x,\lambda)&= \frac{1}{2\pi i} \int_{(\theta)} e^{y\sqrt{x}} \frac{2^{\rho}}{\sqrt{\pi}}\bar{\nu}(T) i^{-2k}  2^{4k}\Gamma\left(2k + \rho + \frac{1}{2}\right)\left(\frac{2\pi}{\lambda} \right)^{2k} \sum^\infty_{m=1} \frac{a_m}{(y^2 + (\frac{4\pi}{\lambda})^{2}m)^{2k + \rho + \frac{1}{2}}} \, dy ,\nonumber 
 \end{align}
 $$ h_{2}(x,\lambda)= \frac{1}{2\pi i} \int_{(\theta)}  e^{y\sqrt{x}} \frac{2^{\rho}}{\sqrt{\pi}} \frac{a_{0}\Gamma(2k+\rho+\frac{1}{2})}{y^{4k+2\rho+1}}\nu(T)e^{\pi ik} \left(\frac{8\pi}{\lambda} \right)^{2k} dy,$$
 $$ h_{3}(x,\lambda)= \frac{1}{2\pi i} \int_{(\theta)} e^{y\sqrt{x}} \frac{2^{\rho+1}}{\sqrt{\pi}} a_{0}\frac{\Gamma(\rho+\frac{1}{2})}{y^{2\rho+1}} dy,$$  and
 $$ h_{4}(x,\lambda)= \frac{1}{2\pi i} \int_{(\theta)} e^{y\sqrt{x}} \frac{2^{\rho}}{\sqrt{\pi}} \frac{\Gamma(2k+\rho+\alpha_{j}+\frac{1}{2})}{y^{2\rho+1}}\left(\frac{8\pi }{\lambda y^{2}} \right)^{2k} \nu(T)e^{\pi ik} \sum_{j=1}^{N} \left(\frac{8\pi i}{\lambda y^{2}} \right)^{\alpha_{j}} \sum_{t=0}^{M_{j}}\beta_{jt} \left(\frac{i\pi}{2} \right)^{t}dy.$$
 The evaluation of $h_{1}(y,\lambda)$~follows as in \cite{CN}. The inter change of integration and summation being justified for $\rho>\beta-2k+\frac{1}{2}$ in the right-hand side of $h_{1}(x,\lambda).$~Applying the integral formula (see \cite{CN}, page $22$)
 \begin{equation}\label{PR8}
 \frac{1}{2\pi i}\int_{w-i\infty}^{w+i\infty}\frac{e^{bs}}{(a^{2}+s^{2})^{\nu+\frac{1}{2}}}ds=J_{\nu}(ab)\left(\frac{b}{a}\right)^{\nu}\frac{\sqrt{\pi}2^{-\nu}}{\Gamma(\nu+\frac{1}{2})}, 
 \end{equation}
 where $w>0,\nu>-\frac{1}{2},a>0,b>0.$ 
  After simplifying each expressions, we obtain
 \begin{equation}\label{h1}
 h_{1}(x,\lambda)=\bar{\nu}(T)e^{-\pi ik}\left(\frac{\lambda}{2\pi} \right)^{\rho}2^{-\rho}\sum_{m=1}^{\infty}a_{m}\left(\frac{x}{m} \right)^{\frac{2k+\rho}{2}} J_{2k+\rho}\left( \frac{4\pi\sqrt{mx}}{\lambda}\right).
 \end{equation} 
 To evaluate $h_{2}(x,\lambda)$ we rewrite as $$h_{2}(x,\lambda)=\frac{2^{\rho}}{\sqrt{\pi}}a_{0}\nu(T)e^{\pi ik} \left(\frac{8\pi}{\lambda} \right)^{2k}\Gamma\left(2k+\rho+\frac{1}{2}\right)\frac{1}{2\pi i} \int_{(\theta )} \frac{e^{y\sqrt{x}}}{y^{4k+2\rho+1}} dy.$$ 
 Put $$g_{1}(y)=\frac{1}{2\pi i} \int_{(\theta)} \frac{e^{y\sqrt{x}}}{y^{4k+2\rho+1}} dy. $$ To compute $g_{1}(y)$~we use the integral representation of reciprocal gamma function Hankel's formula in [\cite{WG}, $page 246$] given by
 \begin{equation}\label{PR9}
 \frac{1}{2\pi i}\int_{\sigma-i\infty}^{\sigma+i\infty}e^{t}t^{-z}dt=\frac{1}{\Gamma(z)}, \ \mbox{ where Re}z>0, \sigma>0.
 \end{equation}
 and applying  the  formula in (\ref{PR3}) for $\frac{1}{\Gamma(4k+2\rho+1)}$. After simplification we obtain 
 \begin{equation}\label{h2}
 h_{2}(x,\lambda)=2^{-\rho}a_{0}\nu(T)e^{\pi ik}\left( \frac{2\pi}{\lambda}\right)^{2k}\frac{x^{2k+\rho}}{\Gamma(2k+\rho+1)}.
 \end{equation}
 We evaluate $h_{3}(x,\lambda)$~applying the techniques analogous used in evaluating $h_{2}(y,\lambda).$ Thus, we conclude that
 \begin{equation}\label{h3}
 h_{3}(x,\lambda)=\frac{2a_{0}}{2^{\rho}}\frac{x^{\rho}}{\Gamma(\rho+1)}.
 \end{equation}
 Next we consider $h_{4}(x,\lambda)$ the interchange of the integration and the finite double sums is justified and we write as;
 \begin{align*}
 h_{4}(x,\lambda)=\frac{2^{\rho}}{\sqrt{\pi}} \nu(T)e^{\pi ik} \left(\frac{8\pi}{\lambda} \right)^{2k} \Gamma\left( 2k+\rho+\alpha_{j}+\frac{1}{2}\right)\sum_{j=1}^{N}\left(\frac{8\pi i}{\lambda}\right)^{\alpha_{j}} \sum_{t=0}^{M_{j}}\beta_{jt} \left(\frac{i\pi}{2}\right)^{t}\\
 \times \frac{1}{2\pi i} \int_{(\theta) } \frac{e^{y\sqrt{x}}}{y^{4k+2\rho+\alpha_{j}+1}} dy.
 \end{align*}
 Now put $$g_{2}(y)=\frac{1}{2\pi i} \int_{(\theta)} \frac{e^{y\sqrt{x}}}{y^{4k+2\rho+\alpha_{j}+1}} dy.$$ Applying (\ref{PR9}) in $g_{2}(y)$ and using  formula~(\ref{PR3})~ for $\frac{1}{\Gamma(4k+2\rho+1\alpha_{j}+1)}.$ ~By substituting the expression derived for $g_{2}(y)$ and simplifying, we obtained
 \begin{equation}\label{h4}
 h_{4}(x,\lambda)=2^{-\rho}\nu(T)e^{\pi ik}\frac{x^{2k+\rho}}{\Gamma(2k+\rho+\alpha_{j}+1)}\left(\frac{2\pi}{\lambda} \right)^{2k}\sum_{j=1}^{N}\left(\frac{2\pi ix}{\lambda}\right)^{\alpha_{j}} \sum_{t=0}^{M_{j}}\beta_{jt}\left(\frac{i\pi}{2}\right)^{t}.
 \end{equation}  
 Therefore. we have  
 \begin{equation}\label{h1234}
 \frac{1}{\Gamma(\rho+1)}{\sum_{m\leq x}}^{\prime }a_{m}(x-m)^{\rho}2^{-\rho}=h_{1}+h_{2}+h_{3}+h_{4}.
 \end{equation}
 By substituting the expressions (\ref{h1}), (\ref{h2}), (\ref{h3}), and~(\ref{h4}) in to $h_{1},$~$h_{2},$~$h_{3},$ and $h_{4}$~respectively  and multiplying both sides of  the  equation (\ref{h1234}) by $2^{\rho}$ the identity in (\ref{Efra2}) holds for $\rho\geq 2\beta-2k+\frac{1}{2}.$~This completes the proof of the converse, since (\ref{Efra2}) is equivalent to (\ref{Func2})as  we shown in Theorem \ref{ThM42}. 
  \section{Conclusion and Discussion}
  
  In this paper we have used techniques of Chandrasekharan and Narasimhan and extended the results of Sister Ann M. Heath to entire automorphic integrals involving rational period functions on discrete Hecke group. Moreover, we have established equivalence of two arithmetical identities with a functional equation associated with automorphic integrals involving log-polynomial-period functions on the Hecke groups. As noted in the introduction part of this paper Sister Ann M. Heath considered the functional equation in the Hawkins and Knopp context and showed its equivalence to two arithmetical identities associated with entire modular cusp integrals involving rational period functions for the full modular group.

 \end{document}